\documentclass{article}
\usepackage{ifpdf}
\ifpdf
   \usepackage[pdftex]{hyperref}
   \usepackage{graphicx}
   \DeclareGraphicsExtensions{.pdf,.jpg}
\else
   \usepackage{graphicx}  
   \usepackage{hyperref}
\fi
\ifpdf

\else

\fi

\usepackage{graphicx,amssymb,amstext,amsthm,url}
\newcommand{\vek}[1]{\mathbf{#1}}
\newcommand{\Reals}{\mathbb{R}}

\newcommand{\infin}{first-order }
\newcommand{\Infin}{First-order }
\newcommand{\ol}{\overline}
 \newcommand{\ovtri}{\widehat{G}}

\newtheorem{theorem}{Theorem}

\newtheorem{corollary}{Corollary}

\newtheorem{conjecture}{Conjecture}

\begin{document}

\title{Geometric Properties of Assur Graphs}
\author{Brigitte Servatius\thanks{ Mathematics Department, 
        Worcester Polytechnic Institute.  bservatius@math.wpi.edu} \and %
        Offer Shai\thanks{Faculty of Engineering, Tel-Aviv University. shai@eng.tau.ac.Il}
        \and Walter Whiteley
        \thanks{Department of Mathematics and Statistics, York University
         Toronto, ON, Canada. whiteley@mathstat.yorku.ca. Work supported
        in part by a grant from NSERC Canada}}

\maketitle

\begin{abstract} In our previous paper, we presented the combinatorial 
theory for minimal isostatic pinned frameworks - Assur graphs - which arise
in the analysis of mechanical linkages.   In this 
paper we further explore the geometric properties of Assur graphs, with a focus on
singular realizations which have static self-stresses. We provide a new geometric
characterization of Assur graphs, based on special singular realizations. 
These singular positions
are then related to dead-end positions in which 
an associated mechanism with an inserted driver will stop or jam.   
\end{abstract}


\section{Introduction} 

In our previous paper~\cite{SSW}, we developed the combinatorial properties of a class of graphs which arise 
naturally in the analysis of minimal one-degree of freedom mechanisms in the plane with one driver, 
with one rigid piece designated at the `ground'.  
We defined an underlying isostatic graph (generically independent and rigid) formed by replacing the driver 
- the Assur graph, named after the Russian mechanical engineer who introduced and began the analysis of this class. 
Every other mechanism,  which is independent (whose degree of freedom increases if we remove any bar) is formed by
composing a partially ordered collection of $k$ such Assur graphs 
(see Figure~\ref{assurscheme}).  The techniques of combinatorial rigidity provided an 
algorithm for decomposing an arbitrary linkage into these Assur components.  

In this paper, we investigate the geometric properties (self-stresses and \infin  motions) 
of such Assur graphs $G$, 
when realized as a framework $G(\vek{p})$ in special or singular positions $\vek{p}$. 
The properties of a full self-stress combined with a ull motion relative to the ground, at selected singular positions
provide an additional, geometric 
necessary and sufficient condition for Assur graphs (see \S3).

Geometric properties of Assur graphs are also important in terms of special positions reached by the
mechanisms, when moving under forces applied through a driver (\S4).  
These position include configurations in which this driver is unable to force a continuing motion, because
of the transmission difficulties, or reach  `dead end' positions in which the mechanism 
will `jam' under the motion of the driver, and continued motion will requre will need to be reversed.  

As with the first paper, our initial motivation was to provide a more complete grounding and 
mathematical precision for some geometric observations and conjectures developed by Offer Shai, 
and presented at the Vienna workshop in April/May 2006.  More generally,
 using the \infin  and static theory of plane frameworks, we want to provide 
a careful mathematical description for the 
properties, observations and operations used by mechanical engineers in their practice.  
We also hope to develop new techniques
and extensions in an ongoing collaboration between mathematicians and engineers.

\section{Preliminaries}
We will summarize some key results from the larger literature on rigidity~\cite{GSS,Wh}, 
and from our first paper~\cite{SSW}.  Throughout this paper, we will assume 
that all frameworks are in the plane and we only consider rigidity in the plane. 

\subsection{Frameworks and the rigidity matrix}

A {\em plane framework} is a graph $G = (V,E)$ together with a
configuration $\vek{p}$ of points for the vertices $V$ in the Euclidean
plane, with pairs of vertices sharing an edge in distinct positions.  Together they are written $G(\vek{p})$.
  A {\em \infin motion} of a framework $G(\vek{p})$ is an assignment of plane vectors $\vek{p}'$ 
to the $n=|V|$ vertices of $G(\vek{p})$ such that, for each edge $(i,j)$ of $G$:
  \begin{equation}\label{InfEqn}
      (\vek{p}_i - \vek{p}_j)\cdot(\vek{p}'_i - \vek{p}'_j) = 0.
  \end{equation}
If the only \infin motions are
{\em trivial}, that is, they arise from \infin translations
or rotations of $\Reals^2$,  then we say that the framework is
{\em infinitesimally rigid } in the plane.  

   Equation~\ref{InfEqn} defines a system of linear equations,
indexed by the edges $(i,j)\in E$, in the variables for the unknown velocities $\vek{p}'_i$
for the framework $G(\vek{p})$.  The 
{\em rigidity matrix} $R_{G}(\vek{p})$  is the real $E$ by $2n$
matrix of this system.  As an example,
we write out coordinates of $\vek{p}$ and of the
rigidity matrix $R_{G}(\vek{p})$, in the case $n=4$ and the complete graph $K_{4}$.
  \[  \vek{p} = (\vek{p}_{1},\vek{p}_{2},\vek{p}_{3},\vek{p}_{4}) =
        (p_{11},p_{12},p_{21},p_{22},p_{31},p_{32},p_{41},p_{42}); \]
  \[ \left[\begin{array}{cccccccc}
   _{p_{11}-p_{21}} & _{p_{12}-p_{22}} & _{p_{21}-p_{11}} & _{p_{22}-p_{12}} &
   _{0}             & _{0}             & _{0}             & _{0}              \\

   _{p_{11}-p_{31}} & _{p_{12}-p_{32}} & _{0}             & _{0}             &
   _{p_{31}-p_{11}} & _{p_{32}-p_{12}} & _{0}             & _{0}              \\

   _{p_{11}-p_{41}} & _{p_{12}-p_{42}} & _{0}             & _{0}             &
   _{0}             & _{0}             & _{p_{41}-p_{11}} & _{p_{42}-p_{12}}  \\

   _{0}             & _{0}             & _{p_{21}-p_{31}} & _{p_{22}-p_{32}} &
   _{p_{31}-p_{21}} & _{p_{32}-p_{22}} & _{0}             & _{0}              \\

   _{0}             & _{0}             & _{p_{21}-p_{41}} & _{p_{22}-p_{42}} &
   _{0}             & _{0}             & _{p_{41}-p_{21}} & _{p_{42}-p_{22}}  \\

   _{0}             & _{0}             & _{0}             & _{0}             &
   _{p_{31}-p_{41}} & _{p_{32}-p_{42}} & _{p_{41}-p_{31}} & _{p_{42}-p_{32}}

   \end{array}\right]    \]

   A framework
$(V,E,\vek{p})$, with at least one edge, is infinitesimally rigid (in dimension $2$) if and
only if the matrix of $R_{G}(\vek{p})$ has rank $2n-3$.
We say that the configuration on $n$ vertices $\vek{p}$ is {\em in generic position} 
if the determinant
of any submatrix of $R_{K_{n}}(\vek{p})$ is zero only if it is identically
equal to zero in the variables $\vek{p}_i$.
For a generic configuration $\vek{p}$, linear dependence of the rows of
$R_{G}(\vek{p})$ is determined by the graph and the rigidity properties of a graph
 are the same for any
generic embedding. A graph $G$ on $n$ vertices  is {\em generically rigid} if the rank
$\rho$ of its rigidity matrix $R_G (\vek{p})$ is $2n-3$, where $R_G (\vek{p})$ is
the submatrix of $R(\vek{p})$ containing all rows corresponding to the edges of $G$,
for a generic configuration $\vek{p}$ for $V$. 

An \infin motion $\vek{p}'$ is a solution to the matrix equation
 $R_{G}(\vek{p})\vek{p}' = 0$, and \infin rigidity is studied through the
 column rank.  We can also analyze the rigidity through 
 the row rank of the rigidity matrix, or through the cokernel: 
 $[\vek{\Lambda}]R_{G}(\vek{p}) = \vek{0}$.  Equivalently,
 these row dependencies $\Lambda$ are assignments
 of scalars $\lambda_{{ij}}$ to the edges such that at each vertex $i$:
 \begin{equation}\label{stress}
 \sum_{\{j|(i,j)\in E\}} \lambda_{ij} (\vek{p}_i - \vek{p}_j) = \vek{0}
  \end{equation}
These dependencies $\Lambda$ are called {\em static self-stresses}, or {\em self-stresses} 
for short, in the 
 language of structural engineering and mathematical rigidity, and the cokernel is the
 vector space
 of self-stresses.  Equation~\ref{stress} is also called the {\em equilibrium condition}, since
 the entries $ \lambda_{ij} (\vek{p}_i - \vek{p}_j)$ can be considered as forces applied
 to the vertex $\vek{p}_{i}$.  Equation~\ref{stress} then states that these forces are in local equilibrium
 at each vertex.  
 Equivalently, a framework $G(\vek{p})$  is {\em independent} if the only self-stress is the
 zero stress, and we see that framework with at least one bar is \infin  rigid if and only if 
 the space of self-stresses has dimension $|E|-(2n-3)$. 

\subsection{Isostatic graphs and rigidity circuits}
 A framework $G(\vek{p})$ is {\em isostatic} if it is infinitesimally rigid and 
 independent.   There is a fundamental characterization of generically isostatic graphs
 (graphs that are isostatic when realized at generic configurations $\vek{p}$).
 
\begin{theorem}[Laman~\cite{Laman}]\label{LamanTheorem}
                A graph $G=(V,E)$ has a realization $\vek{p}$ in the plane as an
                 \infin  rigid, independent
                framework $G(\vek{p})$  if and only if 
                $G$ satisfies Laman's conditions:  $|E| = 2|V|-3$
                \begin{equation}\label{LamansIneq}
                    |F| \leq 2|V(F)| - 3
                    \mbox{~for all $F \subseteq E, F \not = \emptyset $};
                \end{equation}
                Such a graph is also generically rigid in the plane. 
 \end{theorem}
 
Minimally dependent sets, or {\em circuits},  are edge sets
$C$ satisfying
$|C|=2|V(C)|-2$ and every proper non-empty subset of $C$ satisfies inequality (\ref{LamansIneq}).
Note that these circuits, called {\em rigidity circuits}, always have an even number of edges.

\begin{figure}[htb]
\centering
\includegraphics{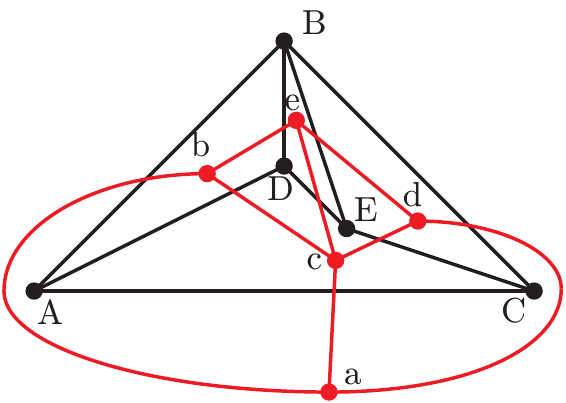}
\caption{The geometric dual $G^{*}$ (red) of a planar graph $G$ (black)\label{recip05Fig}}
\end{figure}

If a rigidity circuit $C=(V,E)$ induces a planar graph, then a planar
embedding of $C$ (with no crossing edges) has as many vertices as it has faces, which follows
immediately from Euler's formula for planar embeddings ($|V|-|E|+|F|=2$). 
The geometric dual graph $C^d$ is also a
rigidity circuit, see~\cite{SC, CMW}. (The vertices of $C^d$ are the
faces of the embedded graph $C$ and two vertices of $C^d$ are
adjacent if the corresponding faces of $C$ share an edge, see
Figure~\ref{recip05Fig}).

In our previous paper~\cite{SSW}, we presented an overview theorem which presented
an inductive construction of all rigidity circuits from $K_{4}$ by a simple sequence
of steps, along with some other properties of circuits.

\subsection{Reciprocal diagrams}
\begin{figure}[htb]
\centering
\includegraphics{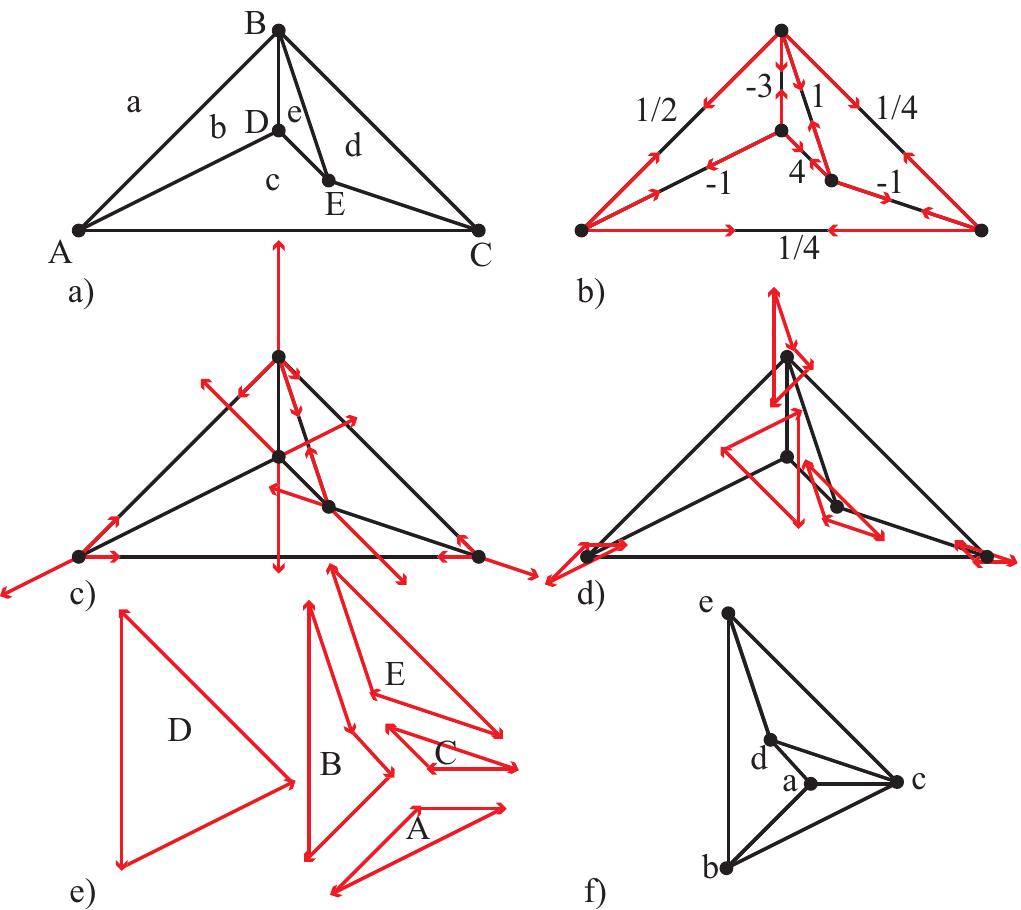}
\caption{The geometric reciprocal (f) of a self-stress on a planar 
framework (c).\label{recip04Fig}}
\end{figure}

In this paper, we will use a classical geometric method for analyzing self-stresses in 
planar frameworks (frameworks on planar graphs): 
the reciprocal diagram~\cite{CW, CW2}.  This construction has 
a rich literature in structural engineering, beginning with the work of 
James Clerk Maxell~\cite{Maxwell} and continuing with the work of 
Cremona~\cite{Cremona}.  
This technique has been revived in the last 25 years as a valuable technique for 
visualizing the behaviour
of such frameworks~\cite{CW, CW2}, in specific geometric realizations, as well as for the study of
mechanisms~\cite{Shai2, Shai3}.   We describe this construction 
and some key properties in the remainder of
this subsection.   An example is worked out in Figure~\ref{recip04Fig}.

Given a framework $G(\vek{p})$ with a non-zero self-stress $\Lambda$ 
(Figure~\ref{recip04Fig} (a)-(c)),
there is a geometric way to verify the vertex equilibrium conditions of 
Equation~\ref{stress}.  If the forces $\lambda_{ij} (\vek{p}_i - \vek{p}_j)$
are drawn end to end, as a polygonal path, then
$$\sum_{\{j|(i,j)\in E\}} \lambda_{ij} (\vek{p}_i - \vek{p}_j) = \vek{0}$$
if, and only if, the path closes to a polygon (classically called the polygon of forces, 
Figure~\ref{recip04Fig} (d), (e)). 
If we start with a self-stress on a planar drawing of a planar graph $G$, 
then we can cycle clockwise through the edges at a vertex in the order of the edges, creating
a polgon for the original vertex, and a vertex for each of the `faces' (regions)
of the drawing.  When we create polygons for each of the original vertices,
we note that the two ends of each original bar use opposite vectors 
(Figure~\ref{recip04Fig}~(b), (c)): 
$\lambda_{ij} (\vek{p}_i - \vek{p}_j)$ and $\lambda_{ij} (\vek{p}_j - \vek{p}_i)$.
We can then patch these polygons together at each dual vertex (original region of the drawing)
(Figure~\ref{recip04Fig}~(e), (f)).
Overall, there is a question whether all of these local polygons patch together 
into a global drawing of the dual graph of the planar drawing?  If we started with a 
self-stress in a planar graph, the answer is yes
 \cite{CW, CW2}!  As just described, this is a {\em reciprocal figure} as studied by 
 Cremona~\cite{Cremona}, with original edges parallel to the edges in the reciprocal drawing.  
  Note that only the edges
 with a non-zero self-stress $\lambda_{ij}$ will have non-zero length in the reciprocal.  
 If the entire reciprocal is turned $90^{\rm o}$, then we have the reciprocal figures as
 presented by Maxwell~\cite{Maxwell, CW, CW2}. 

Conversely, we can start with a configuration for a planar graph $G(\vek{p})$
and a configuration for the dual graph $G^{d}(\vek{q})$ with edges in 
 $G(\vek{p})$ parallel to the edges in $G^{d}(\vek{q})$: a {\em reciprocal pair}.  
 We can then use the 
 lengths of the dual edge to define the scalars for a self-stress, by solving for $\lambda_{ij}$ in
 the equation:
 $$\lambda_{ij} (\vek{p}_i - \vek{p}_j)
 =  (\vek{q}_h - \vek{q}_k)$$
 where edge $hk$ is dual to the edge $ij$. 
 Overall, the existence of the reciprocal also implies the existence of a self-stress and this
 self-stress would recreate the reciprocal.  Note that, in this presentation each of the 
 drawings is a reciprocal of the other.  That is, each side corresponds to a self-stress
 of the other side of the pair. 

The following theorem summarizes some key properties of frameworks and their reciprocals.  
We translate the results of \cite{CW,CW2} into equivalent statements of the form 
we will use in \S3. 

\begin{theorem}[Crapo and Whiteley~\cite{CW2}] Given a 
planar framework $G(\vek p)$ with a self-stress, 
and a reciprocal diagram $G^{d}(\vek{ q})$, there are isomorphisms between:
\begin{enumerate}
\item  the vector space of self stresses of the reciprocal framework $G^{d}(\vek{ q})$;
\item  
the space of frameworks $G(\vek {p}^{||})$ with this reciprocal (with one fixed vertex);  and 
\item 
the space of parallel drawings  $G(\vek {p}^{||})$ with one vertex fixed (equivalently, \infin motions with
one vertex fixed) 
of the framework $G(\vek p)$.
\end{enumerate}
\end{theorem}

Reciprocal diagrams are particularly nice for rigidity circuits, as they exist for all generic realizations.
They were studied
in a much broader context in ~\cite{CW,WW1}  as well as in the context of linkages in~\cite{Shai3}.
We will return to them in the proof of Theorems~\ref{planarselfstress} and~\ref{selfstress}.

\subsection{Isostatic pinned frameworks}
Given a framework, we are interested in its internal motions, not the trivial ones,
so we  `pin' the framework by prescribing, for example, the coordinates of the endpoints of an
edge, or in general by fixing the position of the vertices of some rigid subgraph. Alternatively, we take some rigid
subgraph (a single bar or an isostatic block) which we make into the `ground' and fix all its vertices which 
connect to the rest of the graph.  We call
these vertices with fixed positions {\em pinned }, the others {\em unpinned}, {\em free}, or {\em inner}.
Edges between pinned vertices are irrelevant to the analysis of a pinned framework.  
We will denote a pinned graph by $G(I,P;E)$, where $I$ is the set of inner vertices, $P$ is the
set of pinned vertices, and $E$ is the set of edges, and each edge has at least one endpoint in $I$.

A pinned graph $G(I,P;E)$ is said to satisfy the
{\em Pinned Framework Conditions} if $|E|=2|I|$ and for all subgraphs
$G'(I', P'; E')$ the following conditions hold:
\begin{enumerate}
\item          $|E'|\leq 2|I'|$ if $|P'| \geq 2$,

\item          $|E'|= 2|I'| -1$ if $|P'|=1$ ,  and

\item          $|E'|= 2|I'|-3$ if $P' = \emptyset$.
\end{enumerate}

We call a pinned graph $G(I,P;E)$ {\em (pinned) isostatic} if $\ol{G}(V, E \cup F)$ is
isostatic as an unpinned graph, where $V\supseteqq I\cup P$, no $F$ has any vertex from $I$ as
endpoint and the restriction of $\ol{G}$ to $P=V\setminus I$ is rigid. In other words, we can ``replace" the
pinned vertex set by an isostatic graph containing the pins and call $G(I,P;E)$ isostatic, if this
replacement graph on the pins produces an (unpinned) isostatic graph $\ol{G}$. Which isostatic framework we 
choose, and whether there are additional vertices there, is not relevant to either the combinatorics in the first paper
or the geometry in this second paper.  In \cite{SSW} we proved the following result:

\begin{theorem}\label{PinThm}
Given a pinned graph $G(I,P;E)$, the following are equivalent:

(i) \ \ \ There exists an isostatic realization of $G$.

(ii) \ \ The Pinned Framework Conditions are satisfied.

(iii) \ For all placements $\vek{P}$ of $P$ with at least two distinct locations and all generic
positions of $I$ the resulting pinned framework is isostatic.
\end{theorem}

\subsection {Combinatorial characterizations of Assur graphs} 
In our previous paper~\cite{SSW}, we proved the equivalence 
of a series of combinatorial properties which became the definition of
an {\em Assur graph}.    

\begin {theorem}\label{CharacterThm}
Assume $G= (I,P; E)$ is a pinned isostatic graph.  Then the following are equivalent:

(i) \ \ $G= (I,P; E)$  is minimal as a pinned isostatic graph:  that is for all
proper subsets of vertices  $I'\cup P'$,  $I' \cup P'$ induces a pinned subgraph
$G' = (I'\cup P',  E')$ with  $|E'| \leq  2|I'| -1$.

(ii) \    If the set $P$ is contracted to a single vertex $p^*$, inducing the unpinned
graph $G^*$ with edge set $E$, then $G^*$ is a rigidity circuit.

(iii)   Either the graph has a single inner vertex of degree~$2$ or each time we
delete a vertex, the resulting pinned graph has a motion of all inner vertices
(in generic position).

(iv) \ Deletion of any edge from $G$ results in a pinned graph that has a motion
of all inner vertices
(in generic position).
\end {theorem}

An {\em Assur graph} 
is a pinned graph satisfying one of these four equivalent conditions.
\begin{figure}[htb]
\centering
\includegraphics{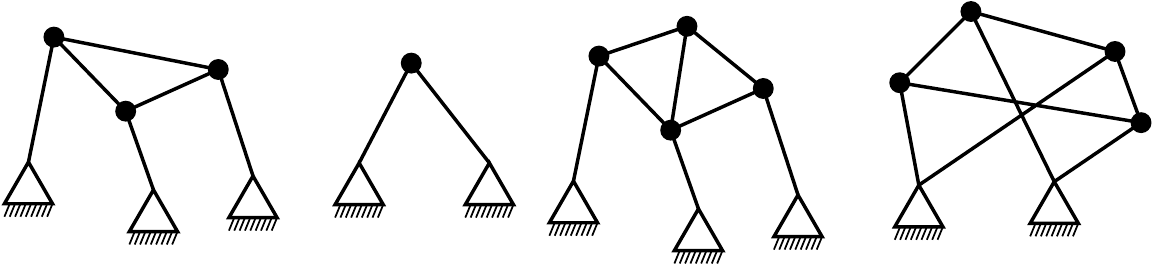}
\caption{Assur graphs\label{assurwordfig01}}
\end{figure}
\begin{figure}[htb]
\centering
\includegraphics{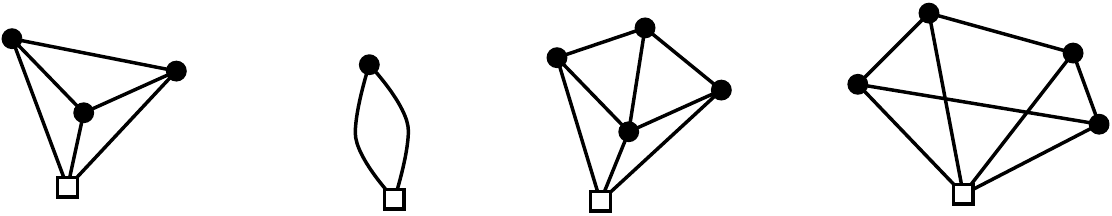}
\caption{Corresponding circuits for Assur graphs\label{assurwordfig03}}
\end{figure}
Some examples of Assur graphs are drawn in Figure~\ref{assurwordfig01} and their corresponding rigidity
circuits in Figure~\ref{assurwordfig03}. 

We also demonstrated a decomposition theorem for all isostatic pinned frameworks.  

\begin{theorem} A pinned graph is isostatic if and only if it decomposes \label{POTheorem}
into Assur components. The decomposition into Assur components is unique.
\end{theorem}

\begin{figure}[htb]
\centering
\includegraphics{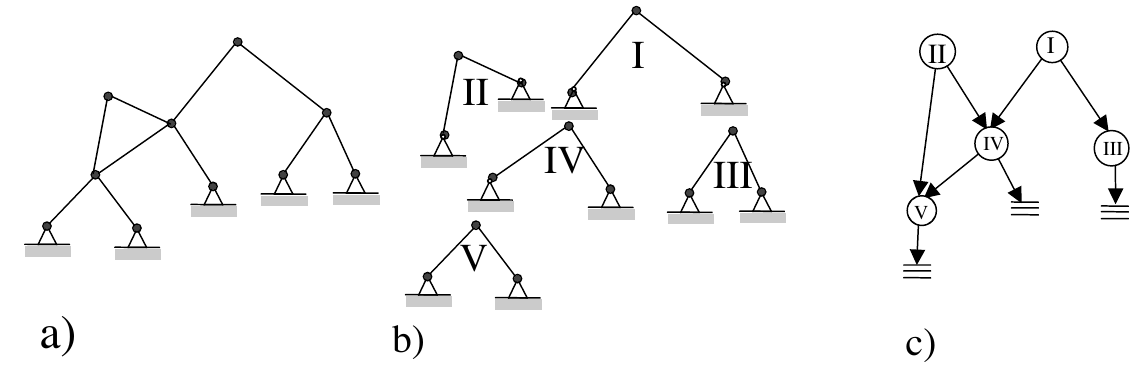}
\caption{An isostatics pinned framework a) has a unique decomposition into Assur graphs b) which is
represented by a partial order or Assur scheme c).\
\label{assurscheme}}
\end{figure}

The decomposition process described in the proof of 
Theorem~\ref{POTheorem} of  \cite{SSW} induces a
partial order on the Assur components of an isostatic graph and this partial order in turn can be used
to re-assemble the graph from its Assur components.  This partial order can be represented in an `Assur 
scheme' as in Figure~\ref{assurscheme}

\section{Singular Realizations of Assur Graphs}

We now show that these Assur graphs have an additional geometric property at selected special positions, 
and that
this geometric property becomes another equivalent characterization.
 From the general analysis of frameworks, we know that the positions $\vek{p}$ of such Assur graphs 
such that $G(\vek{p})$ is not isostatic are the
 solutions to a polynomial {\em Pure Condition}~\cite{WW1}.
 For an pinned isostatic framework, this pure condition is created by
  inserting distinct variables for the coordinates of
 the vertices (including the pinned vertices) and taking the determinant of the
 square $|E|\times2|V|$ rigidity matrix.  
 We are particularly interested in the typical solutions to the pure condition
 (the regular points of the associated algebraic variety).
These properties are related to the behavior of the associated linkage when it reaches
a `dead-end' position and locks ~\cite{Shai2, Penne} (see \S4, Corollary~\ref{deadend} below).

\subsection{A sufficient condition from stresses and motions}
We first show that given a singular realization $G(\vek{p})$ with a special 1-dimensional space of 
self-stresses and 1-dimensional space of \infin motions, $G$ must be an Assur graph.  This is based 
on an observation of Offer Shai.  

\begin{theorem} Assume a pinned graph $G$ has a realization   $p$   such that
\begin{enumerate}
\item  $G(\vek{p})$  has a unique (up to scalar) self-stress $\Lambda$ 
which is non-zero on all edges; and
\item  $G(\vek{p})$   has a unique (up to scalar) \infin  motion, and 
this is non-zero on all inner vertices; 
\end{enumerate}
then  $G$ is an Assur graph.
\end{theorem} 

\begin{proof}  
First we show that $G$ is an isostatic pinned graph.  
For a pinned graph with inner vertices  $V$, and any realization  $\vek{p}$: 
$$|E| - 2|V|  =  \dim ({\text{Stresses}}[G(\vek{p})])  -  
\dim({\text{\Infin \ Motions}} [G(\vek{p})]).$$
Since  $\dim ({\rm Stresses}[G(\vek{p})])  -  
\dim({\text{\Infin \ Motions}} [G(\vek{p})]) = 1-1 =  0$, 
we know that $|E|=2|V|$. 

Similarly, assume there is some subgraph  $G'$
(unpinned, or pinned with one vertex, or pinned with two vertices) which is generically dependent.   
This will always have a non-zero internal self-stress in $G(\vek{p})$ - which is zero outside of this subframework
 $G'(\vek{p})$.  
Therefore,  this unique self-stress cannot be non-zero on all edges.  
This contradiction implies the overall graph $G$ is isostatic.

Now we assume that the graph $G$ is not an Assur graph.  Therefore  
 it can be decomposed into a base Assur graph $G_{A}$, and the rest of the vertices and edges $G_{1}$. 
 With this decomposition, we can sort the vertices and edges of $G_{A}$, 
  to the end of the indicies for the rigidity matrix to give the rigidity matrix
   a block upper diagonal form.  
  Since we have a \infin  motion, non zero on all the inner vertices of $G_{A}$, we have the equation.  
  $$ \Big[\begin{matrix}{ R_{1}(\vek{p}) & R_{1A}(\vek{p})\cr
       0  & R_{A}(\vek{p}) }\end{matrix}\Big]\Big[ \matrix{\vek{u}_{1} \cr \vek{u}_{A}} \Big]= 
       \Big[ \matrix{\vek{0} \cr \vek{0}} \Big] $$
which implies the equation for $G_{A}$:  $[R_{A}(\vek{p})][\vek{u}_{A}]=\vek 0$, 
with $[\vek{u}_{A}]\neq \vek0$ .  
Since $G_{A}$  
is an Assur graph (generically isostatic), $[R_{A}(\vek{p}))]$ has a row dependence.  Equivalently, 
there is a self-stress  $\Lambda_{A} [R_{A}(\vek{p})] = 0$. 
This is also a self-stress on the whole framework $G(\vek{p}))$, 
which is zero on all edges in $G_{1}$. 
Since we assumed $G(\vek{p})$ had a 1-dimensional space of self-stresses, 
this contradicts the assumption that there is a self stress non-zero on all edges.   

We conclude that $G$ is an Assur graph.
\end{proof}

In the next two sections we prove that this condition is also necessary, 
completing this geometric characterization of the Assur graphs.

\subsection{Stressed realizations of planar Assur graphs}
Since this is a geometric theorem, we need to use some key geometric techniques 
for stresses and motions
of frameworks $G(\vek{p})$.   We begin with the special subclass of planar Assur graphs $G$, where we can use 
the techniques of reciprocal diagrams \cite{CW, Shai}.   

\begin{theorem} \label{planarselfstress} If we have a planar Assur graph $G$ then
we have a configuration  $\vek{p}$, such that:
\begin{enumerate}
\item   $G(\vek{p})$ has
a one-dimensional space of self-stresses, and this
self-stress is non-zero on all edges; and 

\item  there is a unique (up to scalar) non-trivial \infin 
motion of  $G(\vek{p})$ and  this is non-zero on all inner vertices.
\end{enumerate}
\end{theorem}

\begin{figure}[htb]
\centering
\includegraphics{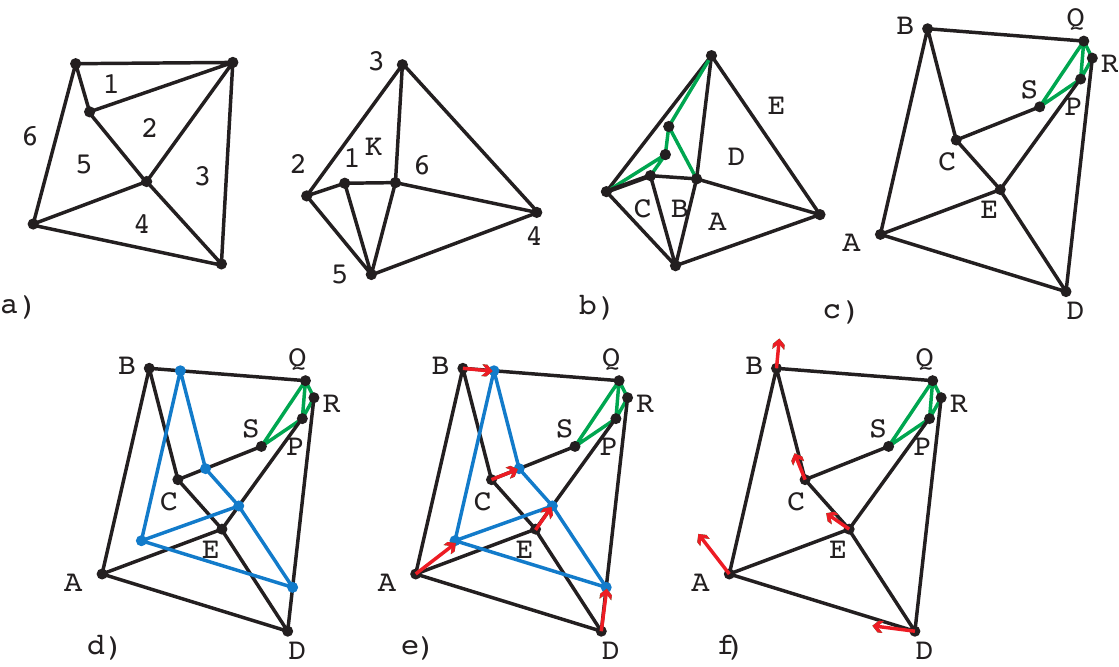}
\caption{The sequence of steps for producing the configuration for a planar Assur graph
which has both a non-zero self-stress and a non-zero motion:   (a) take a generic realization of the
underlying circuit and form its reciprocal; (b) Split the reciprocal face $K$ in order to generate a self-stress
that will separate
the
ground vertex in the original  into predescribed distinct ground vertices (c), still with a self-stress; (d) use a second
self-stress
to form a parallel drawing; (e) use this parallel drawing to create difference vectors; and (f) turn these difference
vectors to create the \infin  motion which is non-zero on all inner vertices.
 \label{StressPlanar}}
\end{figure}

\begin{proof}  
We will use  property (i) from the combinatorial characterization 
Theorem~\ref{CharacterThm}:  $G$ is a 
minimal isostatic pinned framework.  

We assume that the graphs $G$ and $G^{*}$ (with the pinned vertices identified) are planar.
Since $G^{*}$ is a planar circuit, it has a dual graph $G^{*d}$
which is also a planar circuit (Figure~\ref{StressPlanar} (a)).  We take a generic realization of this dual graph 
 $G^{*d}(\vek q^{*})$, which
will have a non-zero self-stress $\Lambda^{*}$ which is non-zero on all edges, and the graph will be \infin  rigid.
We have the corresponding reciprocal diagram $G^{*}(\vek p^{*})$ which also is \infin  rigid and
has a self-stress non-zero on all edges, by the general theory of reciprocal diagrams
(\S2.3 and \cite{CW}).  

\begin{figure}[htb]
\centering
\includegraphics{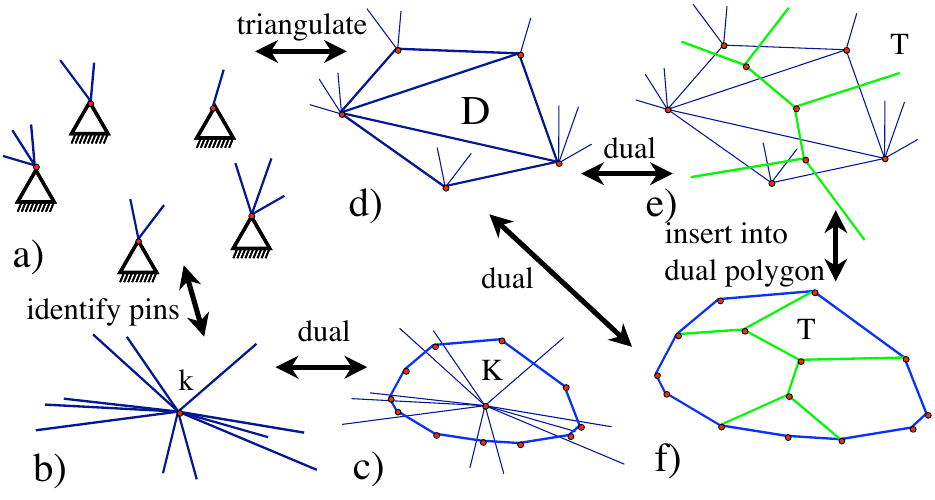}
\caption{Given an Assur graph $G$ we identify the pins to $k$
which has a dual polygon $K$ (a,b,c).  We triangulate the ground (d)
which gives a dual trivalent tree $T$ (e). This tree is inserted into $K$ (f)
giving an additional self-stress whose reciprocal gives back the triangulated 
ground and and separates the pins (f,b).  
 \label{GroundInsert}}
\end{figure}

We will now modify this pair to split up the indentified `ground vertex' $k$ 
while maintaining the self-stress and introducing a \infin  motion $\vek p'$
which is non-zero on all vertices not in the ground.   This process is illustrated in
Figure~\ref{GroundInsert}. 

For simplicity, we create 
this ground for the pinned vertices as an isostatic triangulation on the
pinned vertices.   In the extreme case where we have only two ground vertices, 
  we are adding one edge - and this 
  appears as a corresponding added dual edge $T$ in the reciprocal. 
   More generally, we take the original 
 graph $G$ with $m\geq 3$ pinned vertices, and 
 topologically add an isostatic framework of triangles in place of the `ground' to create 
the extended framework $\ovtri$, with the dual graph  ${\ovtri}^{d}$.   
(For uniformity, we can take a path connecting the 
 pinned vertices $p_{1},\ldots, p_{m}$ and then connect $p_{1}$ to each of the remaining vertices.
 This will be such a generically isostatic triangulation, see Figure~\ref{GroundInsert} (d).) 
 
 If there were $m$ pinned vertices,
 then we add $2m-3$ edges to create the triangulation, and create $t=m-2$ triangles.  
In the dual ${\ovtri}^{d}$, this  
adds a $3$-valent tree $T$ with interior vertices  for each of the triangles 
Figure~\ref{GroundInsert} (e), 
and leaves attached to the vertices of the reciprocal polygon $K$
$K$ (Figure~\ref{GroundInsert} (f)).  Transferring the counts to the reciprocal, 
we have added
$t$ vertices and $2t+1$ edges into the dual polygon $K$.   

Since we added $2t+1$ edges and $t$ vertices to a generically rigid framework 
$G^{*d}(\vek q^{*} )$,
we have added an additional self-stress if all the vertices are in generic position $\vek{q}$.
This added self-stress is non-zero on some of these added edges.  Because the inserted
graph is a $3$-valent tree, if the self-stress is non-zero on one edge, then resolution at 
any interior vertex in general position requires it to be non-zero on all edges at this vertex.
In short, the added self-stress is non-zero on all edges in the tree.

This is now a realization of ${\ovtri}^{d}$ - 
the dual to the original pinned graph with an isostatic triangulated ground 
(Figure~\ref{StressPlanar} (b)). 
In the two dimensional space of self-stresses in the dual, adding a small multiple of the new self-stress
to the original $\Lambda^{*}$ on $G^{*d}(\vek{q}^{*})$ (with zero on the added edges) 
gives a self-stress $\Lambda$ on  
${\ovtri}^{d}(\vek{q})$  
non-zero on all edges.  The reciprocal of 
this self-stress is the desired realization ${\ovtri}(\vek{p})$ of the original pinned framework
with a triangulated (isostatic) ground (Figure~\ref{StressPlanar} (c)).  
Since the self-stress on  ${\ovtri}^{d}(\vek q)$ is non-zero on all edges, all edges are of non-zero
length in ${\ovtri}(\vek{p})$ by the basic properties of reciprocals.   Moreover, since all edges
of ${\ovtri}^{d}(\vek{q})$ have non-zero length, all edges in  ${\ovtri}(\vek{p})$ have non-zero 
self-stress.   With the added sub-framework $D$ replaced
 by the ground, this is the realization $G(\vek{p})$ required for condition (i).  

It remains to prove that this also satisfies condition (ii): there is a non-trivial 
\infin motion with all inner vertices having non-zero velocities while the ground
has zero velocities. 

If we add an additional small multiple of the non-zero self-stress $\Lambda^{*}$  (extended with
zeros on in the added tree in $K$) to $\Lambda$, then we have a second self-stress $\Lambda^{v}$
which is the same on the edges interior to $K$ but different on all other edges.  
Taking a second drawing reciprocal  to $\Lambda^{v}$ will give a second drawing ${\ovtri}(\vek{p}^{||})$
which is identical on the pinned vertices and the ground triangulation, but moves all other 
edges to new positions with different lengths than in ${\ovtri}(\vek{p})$ 
(because of the different self-stress on these edges)
(Figure~\ref{StressPlanar} (d)).
 This is a parallel drawing of ${\ovtri}(\vek{p})$.  In particular, all of the edges of the
 reciprocal polygon $K$ have different stresses, so the edges from inner vertices 
 to the pinned vertices 
 in ${\ovtri}(\vek{p}^{||})$ all have different lengths (Figure~\ref{StressPlanar} (d)).  
 We can take the differences in positions $(\vek{p}^{||} - \vek{p)}$ as 
 parallel drawing vectors $\vek{u}$ (Figure~\ref{StressPlanar} (e)). 
 By general arguments, involving the 90 degree rotation of the `parallel drawing vectors'
  \cite{CW, servatiusWhiteley},  these parallel drawing vectors $u$ convert
 the \infin  motion $\vek{v}=\vek{u}^{\perp}$ of ${\ovtri}(\vek{p})$  
 which is zero on the ground and non-zero on all the inner vertices 
 (Figure~\ref{StressPlanar} (e,f)). 
 This completes the proof
 that $G(\vek{p})$ satisfies condition (ii).  
\end{proof}

\subsection{Extension to non-planar Assur graphs}

We have the desired converse result for all planar Assur graphs.   In order to extend this to singular realizations
of all Assur graphs, we turn to another 19th century technique for converting a singular realization
of a general framework $G(\vek{p})$ into a singular realization of a related planar framework 
without actual crossings of edges, using the same locations
for the original vertices plus the crossing points~\cite{MaxwellBow}.  
This technique is named after the American structural engineer Bow
who introduced it to assist in the analyis of any plane drawing of a framework, in which the 
visible regions and the edges separating them became the pieces for the analysis of the 
framework via a
reciprocal diagram. 

\begin{theorem} \label{selfstress}  If we have an arbitrary Assur graph $G$  then
we have a configuration  $\vek{p}$, such that
\begin{enumerate}
\item   $G(\vek{p})$ has a single self stress, and this
self-stress is non-zero on all edges; and 

\item  there is a unique (up to scalar) non-trivial \infin 
motion of  $G(\vek{p})$ and  this is non-zero on all inner vertices.
\end{enumerate}
\end{theorem}

\begin{figure}[htb]
\centering
\includegraphics{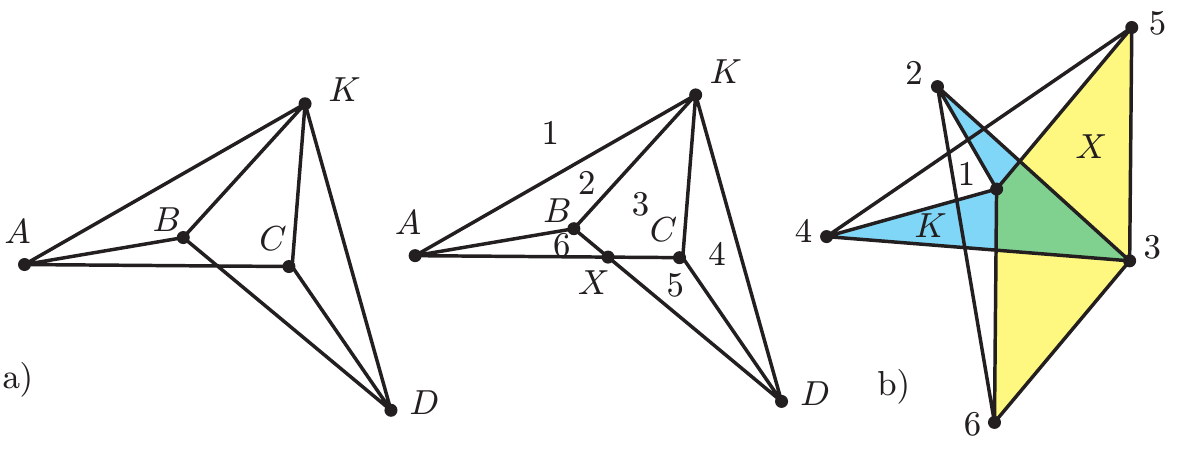}
\caption{Given a non-planar graph (or drawing) (a) we can insert crossing points to create
a planar graph (Bow's Notation).  Working on this planar graph we have a reciprocal (b) which also
is a non-planar drawing of the planar reciprocal (c).
 \label{nonplanarstress1}}
\end{figure}

\begin{proof}  We already know the result for planar Assur graphs.
The key step for non-planar graphs is the classical method called `Bow's notation'.
Given a non-planar framework realizing a graph $G$, we  select pairs of edges with
transversal crossings, and insert those vertices, splitting the two edges, creating a new graph
$G_{b}$~\cite{TW} (Figure~\ref{nonplanarstress1}). 
(Note that these `crossings' do not have to be at internal points of the segments - just not at vertices
of the segments.  The 'crossings are identified topologically, but the added vertices are geometrically
on the points of intersection of the two infinite lines.)  The general
theorem is that the two frameworks have isomorphic spaces of self-stresses, and \infin motions. 

With this technique in mind, we can sketch a plane drawing of the final graph we want,
with the ground triangulation isolated with no crossings.   This sketch identifies the crossing
points to be added, within  the identified circuit - the `Bowed framework'.
Take a generic realization of the identified
circuit $G^{*}$.   Add the crossing points as identified, to create a `planar graph' needed for
the reciprocal diagrams $G_{b}^{*}(\vek{p})$.  Create a reciprocal diagram 
${G_{b}^{*d}}(\vek{q})$
 for this planar framework.  In this reciprocal, 
the identified framework, the duals of the `crossing points' appear as parallelograms.

We now continue with the planar process, as outlined in the previous proof.
With the addition of the vertices and edges to split the  `dual face' $K$ for 
the ground in the dual, we create
a stressed framework, and an extended reciprocal which has the graph of the Bowed framework.
Moreover, since the dual graph is realized with parallelograms dual to the vertices
added in the Bowed framework, the crossings involve transversal crossings with
the required $\times$-appearance for later removal.   
This framework will have a self-stress which is non-zero
on all edges and a non-trivial motion which is non-zero on all inner vertices.  Moreover, this
Bowed framework and the framework with the crossing points removed,  have the isomorphic
spaces of self-stresses and infiniesimal motions.  
In particular, a self-stress which is non-zero on all edges
of the Bowed framework is non-zero on all edges of the original graph,
and the \infin  velocities of the original graph are exactly those velocities assigned 
 at these vertices
within the Bowed framework.

We have created the required configuration for the original (non-planar) graph with a
self-stress non-zero on all edges, and a \infin  motion which is zero on the ground and
non-zero on all free vertices. \end{proof}

\begin{figure}[htb]
\centering
\includegraphics{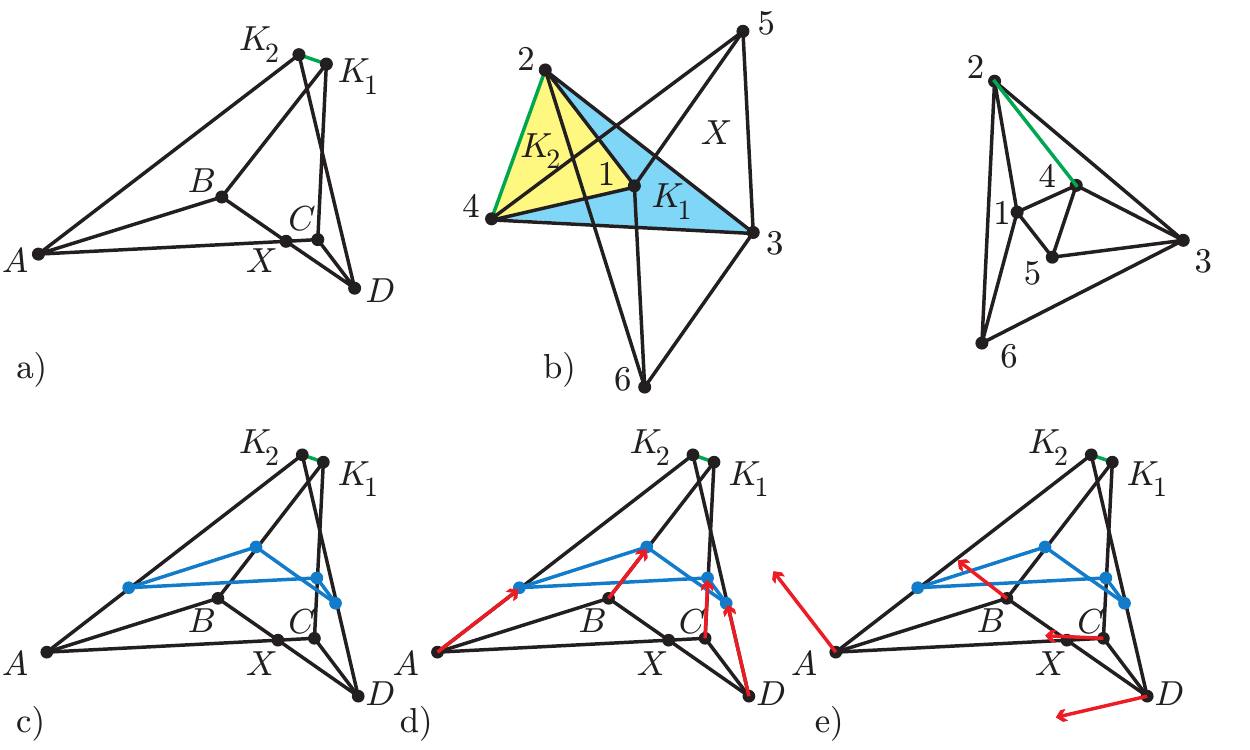}
\caption{Given the reciprocal pair, we can again split the face $K$ (b) and split the ground vertex in
the original.  This configuration of the Bowed graph has a non-zero stress, as does the non-planar
original.   The non-trivial parallel drawing of the Bowed graph is a non-trivial parallel drawing of the
non-planar original (c) and induces the required \infin  motion  on the non-planar original (d,e).
 \label{nonplanarstress2}}
\end{figure}

\subsection{Extensions to other singular realizations}
With some special effort, and careful attention to some geometric details,
it is possible to extend the previous result to show the existence of such 
a special configuration $\vek{p}$ which extends any initial 
configuration of the ground vertices as distinct points.   Without giving all
the details, the idea is to form an isostatic triangulation on the ground vertices as positioned,
which in turn gives an  appropriate dual tree $t$ with dual edges
for the triangles which are on the boundary of the ground, still as rays.  
The `ground' polygon $K$ is then placed on these rays, in general position.
It remains to see
that the rest of the dual graph $G^{*d}$ can be realized with this initial polygon
and a unique self-stress, non-zero on all edges.  Because we started with 
a generic circuit, this can be accomplished by using some details about 
the `polynomial pure conditions' of these graphs~\cite{WW1}, and the occurance
of the remaining vertices in these polynomial conditions.   
\medskip

There is also a conjectured generization of the result above to any realization of an Assur graph
with a non-zero self-stress on all edges.  The key seems to  be not to insist that every
vertex has a non-zero velocity, but to relax this to ask that every bar has
at least one of its vertices with a non-zero velocity.  The reader can review the proof in \S3.1 
to see that this is the condition we actually used in the sufficiency condition.  

\begin{conjecture} \label{selfstressrelax} Assume 
we have an  Assur graph $ G$ and a realization $\vek{p}$
 such that there is a single self stress which is non-zero on all edges.
Then there is a unique (up to scalar) non-trivial \infin 
motion and  this is non-zero on at least one end of each bar.
\end{conjecture}

\section{Inserting Drivers into Assur Graphs}
For simplicitly, we will assume that our graph $G$ is generically independent in the plane.  
In a 1~DOF linkage $G(\vek{p})$ at an independent realization, 
a {\it driver $d$} is either 
\begin{enumerate}
\item [(i)] a {\it piston $ab$} which changes the distance between the pair $ab$
where $ab$ this distance is changing during the 1~DOF motion;
\item[(ii)] an {\it angle driver} which changes the angle $\angle abc$ between two bars 
$ab,bc$ where this angle is changing during during the 1~DOF motion.
\end{enumerate}
More generally, if we have a driver $d$ in aindependent 1~DOF linkage $G(\vek{p})$, 
this driver will cause a finite motion in some  independent realizations, and we continue to call this a piston or 
angle driver 
 even in singular  positions of the same 1~DOF graph as a linkage.  
 We will discuss such singular  positions in \S3.3.
 
\subsection{Replacing drivers} 
In the previous paper~\cite{SSW}, we created the isostatic framework from a 1~DOF  linkage
by `replacing the driver' to remove the degree of freedom.   To return to a 1~DOF
 linkage from an Assur graph, we can `insert a driver'.  
So far, we have used quotation marks here, because we find there are several alternatives for 
the process of replacing the driver, and converse operations to insert a driver because 
the processes are not yet defined - though mechanical engineering practice can
guide us.  We begin be defining one clear process for  `replacing a driver' with an 
added bar.   

\begin{figure}[htb]
\centering
\includegraphics{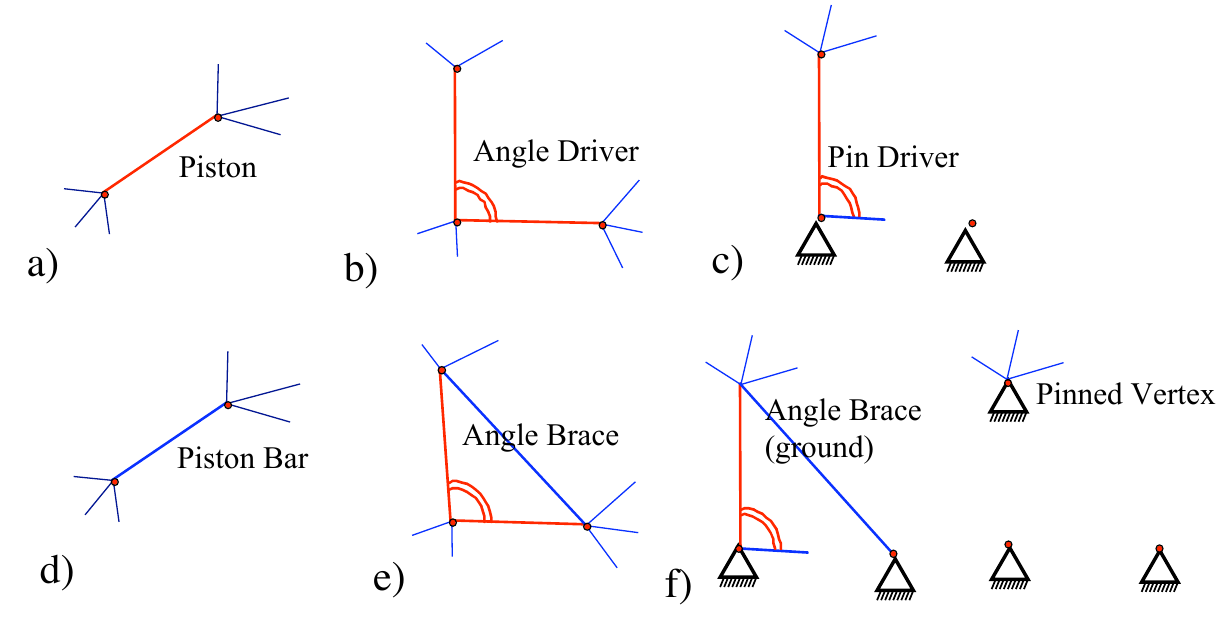} 
\caption{There are four types of drivers: a) driving a distance $ab$ with a piston; 
b) driving an interior angle $\angle abc$; c)  driving an angle at a pin
$\angle a p_{i}p_{j}$, and the special type illustrated in Figure~\ref{ContractInsert}. Below these are the
ways in which each of these drivers is replaced to create an isostatic 
graph  by d) an added edge (for a piston); 
e) an added angle brace; or f) an 
added pinned vertex resulting from adding an angle brace to the ground. 
 \label{driversreplace}}
\end{figure}

A simple method to remove the 1~DOF is to insert one bar in a way that blocks the single 
degree of freedom (Figure~\ref{driversreplace}).  
Specifically, given a one-degree of freedom linkage 
$G(\vek{p})$  at a generic configuration:
\begin{enumerate}
\item  to replace a piston $ab$, we insert the bar $ab$;
\item  to replace an angle driver  on $\angle abc$, at an internal vertex $b$ of 
degree $\geq 3$,  we insert  the bar $ac$;
\item  to replace an angle driver   on $\angle abc$, at an internal vertex $b$ of 
degree $2$,  we insert the bar $ac$ and remove  vertex $b$;
\item  to replace an angle driver on an angle $\angle ap_{i}p_{j}$ where $p_{i}p_{j}$
are pinned vertices,  we add $a$ as a pinned vertex.  
\end{enumerate}
This {\it driver replacement} creates a pinned graph $\ol{G}$ and a pinned framework
 $\ol{G}(\vek{p})$.  

\begin{figure}[htb]
\centering
\includegraphics{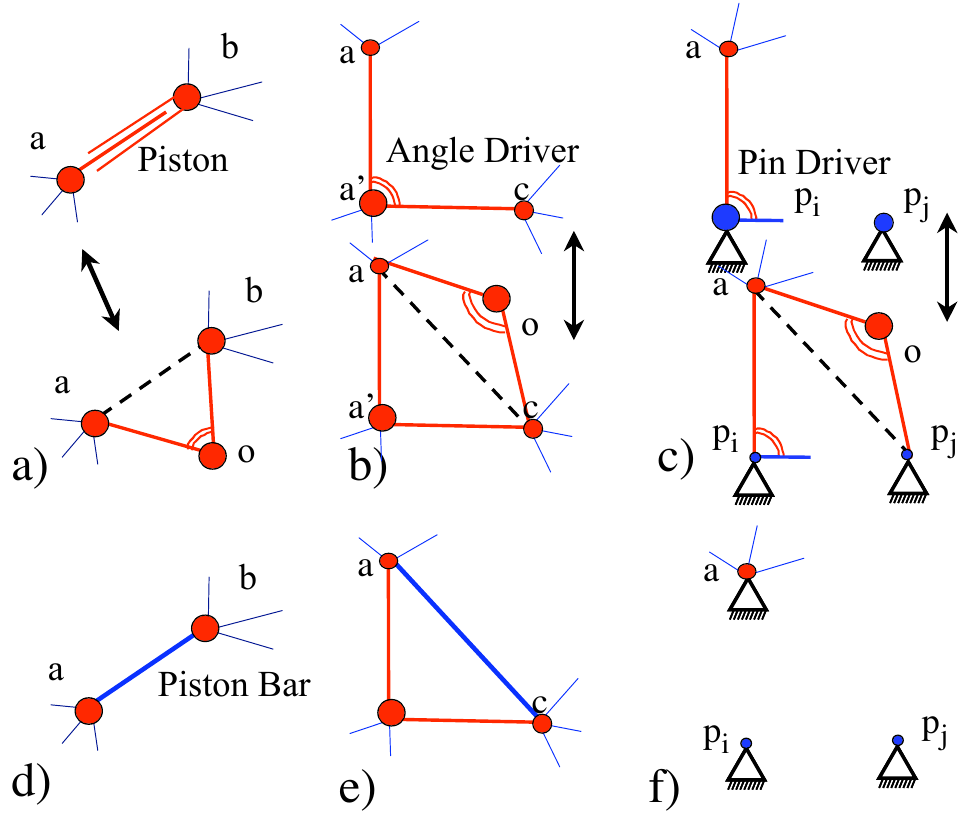} 
\caption{As an alterative to the simple insertions, we can replace each
driver by a 2-valent vertex $o $ (a,b,c), and then contract $o$ to 
one of the end vertices, creating the same result (d,e,f) as the previous insertion. 
 \label{ContractInsert}}
\end{figure}

In the key example of our previous paper~\cite{SSW} Figures~1,2, 
we actually replaced the pistons in two steps:
\begin{enumerate}
\item we replaced the piston  with a $2$-valent vertex attached to the ends
(which is mechanically equivalent); and
\item we contracted one of these edges to 
form the single edge which was inserted above.   
\end{enumerate}
A similar process could be applied to any angle driver, and the net effect would be the
edge insertion presented above (Figure~\ref{ContractInsert}). 

 The driver is {\em active} at a specific position $G(\vek{p})$ if it is possible, infinitesimally, 
 to change the  length of the bar we are adding without changing the lengths of
 any of the other bars.  Specifically, in the rigidity matrix of $\ol{G}$, with the 
 added bar $d$ at the bottom: 
  $R_{\ol{G}}(\vek{p})\vek{p}' =(0,\dots, 0, s_{d})^{tr}$ has a solution $\vek{p}'$
  for all possible values of the strain $s_{d}$ (instantaneous change in length) of  $d$.  
 
More generally, we claim  $\ol{G}(\vek{p})$ is isostatic, and  
$\ol{G}$ is generically isostatic, provided that $G(\vek{p})$ is independent
and the driver was active.  

\begin{theorem}  Given an independent 1~DOF linkage $G$ with an active driver $d$,
the  driver replacement $\ol{G}$ is an 
isostatic pinned framework.   
\end{theorem}

\begin{proof}  Consider an independent realization $G(\vek{q})$.  
There is a 1 dimensional vector space of non-trivial \infin  motions $\vek v'$, with the
pinned vertices fixed (which extends to a finite motion
by general principles of algebraic geometry \cite{Roth}).

For a piston $ab$, the added bar $i,j$ is independent if, and only if,  
the added bar is on a pair $i,j$ with a non-zero {\it strain}:
$$ (\vek q_{i}-\vek q_{j})\cdot (\vek v_{i}- \vek v_{j})\neq 0 $$
The definitions of a piston driver the added bar
has this required property, so inserting this bar blocks the motion, at first-order. 

Similarly, for an angle driver $\angle abc$ at an interior vertex $b$, 
adding the bar $ac$ is also independent, generating an isostatic framework. 

If we were replacing an angle driver at a 2-valent vertex, then with the added bar this vertex is attached
to an isostatistic subframework with just two non-collinear bars.  This vertex can be removed
to leave an isostatic framework on the remaining vertices.   This is done to prevent the appearance
of an extra `Assur component' in the derived isostatic graph and focus the analysis on the behaviour
of the rest of the graph.  

If we were replacing an angle driver $\angle ap_{i}P_{j}$ at the ground, then inserting the bar $a p_{j}$
will create an isostatic framework.  It will also pin the vertex $a$ to the ground, artifically creating
an extra Assur component.  To assist the analysis of the original mechanism, we just pin the vertex
$a$  and analyze the modified pinned framework. 
\end{proof}

We can speak of a 1~DOF graph $G$ with a driver $d$ as an {\it Assur mechanism}, 
if replacing the driver creates an Assur graph $\ol G$.

\subsection{Inserting a Driver}
 Conversely, we can start with an Assur graph, and {\it insert a driver} using one
 of the three steps:
 \begin {enumerate}
 \item[(i))]  Remove a bar $ab$ and insert a piston $ab$;
 \item[(ii)]  Remove an edge $ac$ which is in a triangle $abc$ and insert an angle driver
 on the angle $\angle abc$;
 \item[(iv)]  remove an edge $ac$ and insert a new 2-valent vertex $b$, with bars $ab,bc$
 and an angle driver on the angle $\angle abc$;
 
 \item[(iv)]  if there are at least three pinned vertices $p_{i}p_{j}p_{k}$, make a pinned vertex $p_{k}$
 into an inner vertex $a$, with a single bar to one of the other pinned vertices $p_{i}$
 and an angle driver on the angle $\angle ap_{i}p_{j}$.
 \end{enumerate}
 These operations are the reverse operations of the four ways of replacing a driver.  
 
 As a generic operation, we know that this driver insertion  takes an isostatic framework
 to a 1~DOF framework.  We now show that for an Assur graph, the driver insertion will create
 a 1~DOF framework with all inner vertices in motion
relative to the ground.

\begin{theorem} {\label{driverinsert}} If we have an Assur graph $G$, realized as 
an independent framework 
$G(\vek p)$ with the pins not all collinear,
and we insert a driver as above, then the framework has 1~DOF, with all inner vertices
in motion, and activating the driver will extend this to a continuous path. 
\end{theorem}

\begin{proof}  The original framework is an independent pinned framework with $|E|=2|V|$.
For insertions (i) (ii), we have removed one bar, so there is 1~DOF.  Since this is an Assur graph, 
the Characterization Theorem~\ref{CharacterThm}~(iii) guarantees all inner vertices have a non-zero velocity. 

Moreover, the non-trivial \infin  motion will have a non-zero strain on the pair of the 
removed bar.  If we inserted a piston, this will change this length and drive the motion.  
If we inserted an angle driver on a triangle $abc$, then driving the angle will change the length $ac$
and thus drive the motion. 

Finally, if we changed a pin to an interior vertex, then we can assume that the vertex being made 
an inner vertex is not collinear with two of the other pinned vertices.  We can assume 
that is vertex is 2-valent in the isostatic ground framework, so that making it inner leaves
an isostatic ground, and creates a new inner vertex $a$ attached to the ground by two edges $ap_{i}, ap_{j}$. 
 Removing $ap_{j}$ gives a 1~DOF linkage, as required.  In the resulting motion, $a$ will have
 a non-zero velocity, so driving the angle $\angle ap_{i}p_{j}$ will drive this 1~DOF.  
 It remains to check that all other inner vertices are have non-zero velocities.  
 If some inner vertex $h$ has the zero velocity, then it is attached to the remaining ground through
 an isostatic subframework.  This would mean $h$ is contained in an isostatic pinned subframework
 which does not include $a$.   This is a contradiction of our assumption that $G$ was Assur. 
 \end{proof}

Driver insertion and driver replacement are inverses of one-another.  
That is, if we start with an Assur graph and insert
a driver, then replacing the driver will return us to the same Assur graph.  Conversely, if
we start with a driver, and replace it, then we can choose to re-insert the same driver and
return to the same 1~DOF linkage. (There is a choice in the insertion, one of which is corresponds
to the original replacement. )

However, it is now clear that we have:

(i) \ \ as many ways to insert a piston as we have interior bars;

(ii) \ three time as many ways to insert an angle driver as we have interior triangles;

(iii) as many ways to insert a 2-valent angle driver as there are interior bars;

(iv) \ as many ways to insert a pinned angle driver as we have pins,
 provided there are at least two 
pins.   

\noindent In short, there are a lot of 1~DOF linkage with drivers which come from the same
underlying Assur graph.  All of these will be Assur mechanisms.

\begin{figure}[htb]
\centering
\includegraphics{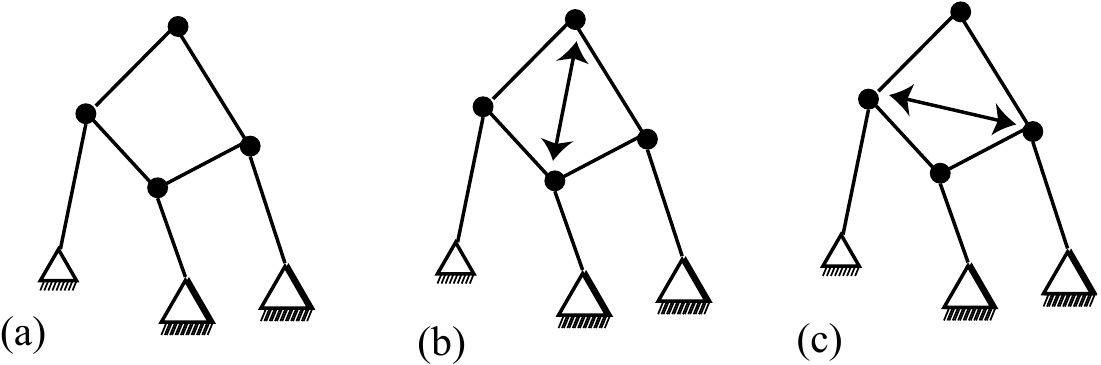} 
\caption{Given a 1~DOF linkage (a), there are several ways to 
insert a driver (b,c).  One of these generates an Assur graph (b), while the
other is a composite (though isostatic) graph (c). 
 \label{DriverInsertAmb}}
\end{figure}

Conversely, if we have a 1~DOF linkage, we can identify a number of 
pairs whose distances are 
change, and angles which are changing.  Each of these could be used to insert a driver.  
However, different insertions will lead, after replacement, to different graphs $\ol{G}$.  
One may be an Assur graph while another may not.    
Figure~\ref{DriverInsertAmb} gives such an example.

\subsection{Singular positions with a driver}

The geometric
 question is: can we find positions at which this still has a self-stress?  If we do, is it possible 
 that some vertices must have zero velocities relative to the ground? Or more 
 generally is it possible that the \infin  motion does not continue in the same direction
 as the original motion?    
 
We recall that if \infin motion $\vek{p}'$ is a \infin motion, then  $-\vek{p}'$  
is also a \infin motion.  If both of 
these velocities extend to a finite motion, then we say the driver has a {\em finite motion 
in both directions}.  
As a contrapositive of the Theorem~\ref{driverinsert}, we have the following corollary.

\begin{corollary}  \label{deadend} If the linkage $G_{d}(\vek{p})$ with the driver does 
not have a finite motion in both directions,  then the linkage with the driver insert
 has a self-stress.  
\end {corollary} 

Such configurations without a finite motion (continuing in both directions) are called 
`dead ends' in the literature of linkages \cite{Shai2}.    For example, the existence 
of a dead-end with an angle driver at the ground requires a self-stress 
with the driver joint pinned, so the original isostatic graph was realized in a singular position.
As another example, if the driver is a piston, 
the driver edge is part of this singular position, so that its
line applies the 'ground force' required for the self-stress of the isostatic graph.  
That an independent 1~DOF linkage can move under a driver to such a singular 
position (with the driver replacement), is one reason why we have investigated 
the occurance of such singular positions with a self-stress on all members (including
the edge used to replace the driver).  That all the inner vertices have non-zero velocities
at that singular position indicates that we could have moved into the self-stress position in 
the driven motion.   

It is not true that every 
self-stress gives a dead end position, just that dead end positions require the self-stress.
The study of such configurations is the subject of a recent paper 
of Rudi Penne \cite{Penne}.  That study focuses on centers of motion, rather
than self-stresses, and it is well-understood in the literature that these are equivalent 
tools for many purposed, each giving its own insight into the geometry and the combinatorics
of linkages.

\section{Concluding comments}
\noindent{\bf Working with several drivers.}
In the previous paper, we presented a decomposition of a general
isostatic pinned framework into Assur components.  With such a decomposition, 
we could insert one driver into all, or some, of the components, creating a larger mechanism 
with as many degrees of freedom as the drivers inserted. See, for example, the mechanism
in Figures~1 and 2 of \cite{SSW}.  These drivers will be independent,
in the sense that each of them could be given distinct instantaneous driving instructions
without any interference or instability.  

We could also insert several drivers into a single Assur graph 
extra analysis will be needed to ensure that these are independent. More generally, given a
mechanism with a number of drivers, their `independence' is equivalent to whether replacing 
all the drivers produces a graph which is isostatic.   

\medskip
\noindent{\bf Projective geometry for self-stresses.}
The instantaneous kinematics and statics of plane frameworks are projectively invariant.  
Thus the singular position of a graph $G(\vek{p})$ can be transferred to any projective image of the configuration
$p$.  In particular, we have seen that it is common in mechanical engineering to include 
pistons (also called `slide joints' in structural engineering).  
These pistons are actually mechanically equivalent to 'joints at infinity' 
between the two ends of the slide \cite{CW} 
- and therefore are incorporated in the geometric (and combinatorial) 
theory we have described in these two papers. 

We also note that spherical mechanisms (with joints built as pines pointed to the 
center of the sphere), share the same projective geometry as the plane mechanisms. 
As this suggests, all of the combinatorial and geometric methods and results 
presented in our two
papers extend immediately from the plane to the spherical frameworks.   

In animating a one degree of freedom mechanism in computer science, it is common to find that
a single link cannot be taken as the `driver' which completes a circle while preserving all of the edges
of the mechanism.   Somewhere along the path, the linkage will experience a self-stress.  Some of these
singular positions will `jam', others will not.  The analysis of the singular positions,  helps clarify this situation.  However the decision 
of which `new driver' to pick to move points along the subject of another study. 

\medskip
\noindent{\bf Extensions to 3-D.} 
In the conclusion of our earlier paper~\cite{SSW}, we indicated that the combinatorial results 
do have appropriate generalizations in 3-space, at least 
for some `nice classes' of structures, such as bar and body structures.    

These will still have drivers and will have special geometry for their dependencies and for 
dead-end positions.   
The specific geometric theorems given are conjectured to 
also extend, but we will need new methods, because techniques such as `reciprocal diagrams' 
are  limited to plane structures.   

\bibliographystyle{plain}
\bibliography{assur1}

\end{document}